\documentclass[english]{amsart}
\usepackage{amssymb,amsmath}
\usepackage{amsfonts,cite,amsthm}
\renewcommand{\le}{\leqslant}
\renewcommand{\ge}{\geqslant}

\newenvironment{Proof}
{\noindent{\bf Proof.}}
{\hfill$\scriptstyle\blacksquare$}

\newenvironment{ProofTh1}
{\noindent{\bf Proof of Theorem~\ref{th1}.}}
{\hfill$\scriptstyle\blacksquare$}

\newenvironment{ProofTh2}
{\noindent{\bf Proof of Theorem~\ref{th2}.}}
{\hfill$\scriptstyle\blacksquare$}

\newtheorem{theorem}{Theorem}
\newtheorem*{col}{Collolary}
\newtheorem{lemma}{Lemma}

\theoremstyle{definition}

\title[Logarithmic Derivatives of Least Deviation from Zero]{Logarithmic Derivatives\\ of Least Deviation from Zero}
\author{Petr Chunaev}
\address{Centre de Recerca Matem\`{a}tica, Campus de Bellaterra, Edifici C, 08193 Bellaterra (Barcelona) Spain}
\email{chunayev@mail.ru}

\begin{document}
\begin{abstract}
We study least deviation of logarithmic derivatives of real-valued polynomials with a fixed root from zero
on the segment $[-1;1]$ in the uniform norm with the weight $\sqrt{1-x^2}$ and without it. Basing on results of Komarov and Novak and on a certain determinant identity due to Borchardt, we also establish a criterion for best uniform approximation of continuous real-valued functions by logarithmic derivatives in terms of a Chebyshev alternance.
\end{abstract}

\maketitle

\section{Introduction}

\textit{Logarithmic derivatives  of algebraic polynomials} (abbreviated to \textit{l.d.}s in what follows) are rational functions of the form
\begin{equation}
\label{n.d.}
\rho_0(z)\equiv 0,\qquad \rho_n(z)=\frac{P'_n(z)}{P_n(z)}=\sum_{k=1}^n\frac{1}{z-z_k}, \qquad z,z_k \in \mathbb{C},\qquad n\in \mathbb{N},
\end{equation}
where $P_n(z)=\prod_{k=1}^n(z-z_k)$ and $n$ is the \textit{degree} of the \textit{l.d}.
Note that \textit{l.d.}s are also widely called \textit{simple partial fractions}.

In rational approximation theory, attention to \textit{l.d.}s was paid by works of Macintyre and Fuchs~\cite{Macintyre}, Gonchar~\cite{Gonchar} and Dolzhenko~\cite{Dolzhenko}.
Later on,  approximating properties of \textit{l.d.}s were studied by Korevaar~\cite{Korevaar} and Chui~\cite{Chui}, who
proposed a construction of \textit{l.d.}s for approximation of analytic functions belonging to the Bergman-Bers class on simply connected domains. Besides, the poles of the \textit{l.d.}s were chosen on the boundaries of
the domains. Such an approximation was motivated by the fact that \textit{l.d.}s specify plane electrostatic fields, consequently, it can be considered as a determination of locations $z_k$ of electrons creating a given field to a high accuracy.

Further investigations of approximating properties of \textit{l.d.}s were initiated by Gorin~\cite{Gorin}, who brought up a problem of least deviation of \textit{l.d.}s with a fixed pole
from zero on the real line. At various times this was studied  by many mathematicians (see \cite{Dan1994} and the references  given there).
Then the following analogue of classical Mergelyan's theorem on polynomial approximations was proved for \textit{l.d.}s with free poles: every function continuous in the compact set $K\subset{\mathbb C}$ with connected complement and analytic in the interior of $K$ can be approximated uniformly
on $K$ with \textit{l.d.}s \cite{DanDan2001}. In many respects this result entailed current intensive study of approximating properties of \textit{l.d.}s and several their modifications on bounded and unbounded subsets of the complex plane (see  \cite{Dan2006,DanKond2010,Kosuh,KomarovNeed,DanChu2011,DDan2001,KomarovPMA2013,Kom2013,Komarov2014,Novak,Chunaev2010,Borodin,Kosukhin,Protasov,Kayumov}).

Hereinafter we consider only real-valued \textit{l.d.}s $\rho_n=\rho_n(x)$, $x\in \mathbb{R}$ (sets of their poles
are symmetric with respect to the real line). Moreover, from now on we set
$$
\|\rho_n\|:=\max_{x\in I}|\rho_n(x)|, \qquad \|\rho_n\|^*:=\max_{x\in I}|\sqrt{1-x^2}\rho_n(x)|, \qquad I:=[-1;1].
$$
Main results (\S\S\ref{s1}--\ref{s-Col}) of our paper are devoted to the following problem stated in \cite{Dan2006}.

\medskip

{\it
Among all real-valued \textit{l.d.}s of degree less than or equal to $n$ with a fixed pole find the \textit{l.d.}s $\widetilde{\rho}_n$ and $\widetilde{\rho}^*_n$ of least deviation from zero on $I$ in the norms $\|\cdot\|$ and $\|\cdot\|^*$.
}

\medskip

This extremal problem is one of possible analogues of classical Chebyshev's question about unitary polynomials of least deviation from zero on $I$ and an analogue of Gorin's problem stated above in the case of a finite segment. We solve it in the class, being denoted by $\mathfrak{R}_n(a)$, of (real-valued) \textit{l.d.}s with the fixed real pole $x=a$.
Note that all results for $I$ can be extended on the case of any segment with the substitution $y=\mu x+\nu$, $\mu>0$. At that it should be taken into account that $\rho_n(x)=\mu \rho_n(y)$, therefore $\|\rho_n\|_{I_{\mu,\nu}}=\mu^{-1}\|\rho_n\|_{I}$, where $I_{\mu,\nu}:=[-\mu+\nu;\mu+\nu]$.

Basing on results of Komarov \cite{Komarov2014} and Novak \cite{Novak} and on a certain determinant identity due to Borchardt \cite{Borchardt}, we also establish a criterion for best uniform approximation of continuous real-valued functions by \textit{l.d.}s in terms of a Chebyshev alternance\footnote{
Let us recall several definitions. Let $f=f(x)$ be a continuous function in $I$.
\textit{Least deviation} of \textit{l.d.}s $\rho_n$ from $f$ on $I$ is the quantity
$\inf\nolimits_{\rho_n}\|f-\rho_n\|$, where the infimum is taken over all \textit{l.d.}s $\rho_n=\rho_n(x)$. The \textit{l.d.}  $\widetilde{\rho}_n$, least deviating from $f$ on $I$, is called a \textit{l.d. of best uniform approximation of $f$ on~$I$}. It is said that $n$ pairwise distinct points $t_k$ in $I$ form \textit{a $($Chebyshev$)$ alternance} of the difference $f-\rho_n$ on $I$ if $f(t_k)-\rho_n(t_k)=\pm (-1)^k \|f-\rho_n\|$, $k=\overline{1,n}$.
} (see \S\ref{s3} for more details). %Note that earlier similar problems concerning connection between an alternance and best approximation were considered by  Ya. Novak \cite{Novak}, but in a much more narrow class of \textit{l.d.}s with real poles.

\section{\textit{L.d.}s of least deviation from zero on $I$ in the norm $\|\cdot\|^*$ }
\label{s1}
In \cite{Dan2006}, there was constructed the \textit{l.d.}
 \begin{equation}
\label{n.d.Dan}
\varsigma_n(A;x)=\frac{2nw}{w^2-1}\frac{A(w^{2n}-1)}{(Aw^n-1)(w^n-A)}, \quad x=\frac{1}{2}\left(w+\frac{1}{w}\right),\quad
w=e^{i\varphi},
%=-\frac{2nA}{\sin \phi}\frac{\sin n \phi}{1-2A\cos n\phi +A^2},
\end{equation}
where $A>1$, $x\in I$, $\varphi \in \mathbb{R}$. It was also shown there that the function $\sqrt{1-x^2}\varsigma_n(A;x)$ has
an alternance consisting of $n$ point in $I$, and $\|\varsigma\|^*={2nA(A^2-1)^{-1}}$. The presence of the alternance suggests that the \textit{l.d.} (\ref{n.d.Dan}) is actually the one, which is a solution to the problem stated above with the norm $\|\cdot\|^*$. And so we prove in this section.
We use the following important property of the \textit{l.d.}~(\ref{n.d.Dan}) below: all its poles belong to the ellipse
\begin{equation}
\label{ellipseChu}
E_p:=\{z:z=p\;\cos t + \sqrt{p^2-1}\;\sin t\;i, \quad t\in [0;2\pi), \quad p>1\},
\end{equation}
where $p=p(A):=\frac{1}{2}\left(A^{1/n}+A^{-1/n}\right)$ is a parameter. This can be  easily deduced from the denominator of (\ref{n.d.Dan}).

As usually, $T_n(x)=\cos n\arccos x$ stand for the Chebyshev polynomials of the first kind (with a leading coefficient $2^{n-1}$), and  $U_n(x)=\frac{1}{n+1}T_{n+1}'(x)$ do for the Chebyshev polynomials of the second kind (with a leading coefficient $2^n$). Properties of the polynomials are well studied (see \cite{Mason}), therefore  we take into account below without any additional explanations for example  that $|T_n(x)|\le 1$ for $x\in I$, ${T_n(1)=1}$, $T_n(-1)=(-1)^n$, and $|T_n(x)|>1$ for $x\in\mathbb{R}\setminus I$.

We now formulate the main result of the section.
\begin{theorem}
\label{th1}
For $n\ge 1$,  the \textit{l.d.}
\begin{equation}
\label{n.d.*}
\tilde{\rho}^*_{n}(a;x)=\frac{U_{n-1}(x)}{\int_a^xU_{n-1}(t)dt}=\frac{nU_{n-1}(x)}{T_{n}(x)-T_n(a)}, \qquad a>\sqrt{2}, %\qquad(a>1),
\end{equation}
is least deviating from zero on $I$ in the norm  $\|\cdot\|^*$ among all those belonging to $\mathfrak{R}_n(a)$. Moreover,
$$
\|\tilde{\rho}^*_{n}\|^*=\frac{n}{\sqrt{T_n^2(a)-1}}.
$$
\end{theorem}

It is easy to check that the \textit{l.d.}s (\ref{n.d.Dan}) and (\ref{n.d.*}) coincide with each other. Nevertheless, we use our representation (\ref{n.d.*}) for convenience of proofs.

Let us first prove several auxiliary results.

\begin{lemma}
\label{lemma4}
For $a>1$ and $n\ge 1$, the function $R_n(x):=\sqrt{1-x^2}\tilde{\rho}^*_{n}(a;x)$ has the following properties.

$(a)$ For $x_k=\cos \frac{\pi k}{n}$, $k=\overline{0,n}$, we have $R_n(x_k)=0$.

$(b)$ The points $a_k$, $k=\overline{1,n}$, being roots of the equation
       \begin{equation}
       \label{equ_(b)}
       T_n(x)=1/T_n(a).
       \end{equation}
form an alternance on $I$ of the function $R_n$, and
       $$
       |R_n(a_k)|=\frac{n}{\sqrt{T^2_n(a)-1}}, \qquad k=\overline{1,n}.
       $$

$(c)$ The poles of $R_n$ belong to the ellipse $E_a$ of the form $(\ref{ellipseChu})$.
\end{lemma}
\begin{Proof} It is clear that the numerator of $R_n$ vanishes both in the endpoints of $I$ and in the roots of the polynomial $U_{n-1}$,
\begin{equation}
\label{Cheb1}
 U_{n-1}(x)=0 \quad \Leftrightarrow \quad x=x_k=\cos \tfrac{\pi k}{n}, \qquad k=\overline{1,n-1}.
\end{equation}
This is $(a)$.

To prove $(b)$ we first find the extreme points $e_k$ of $R_n$ in the interval $I_o:=(-1;1)$. By well-known properties
\begin{equation}
\label{Cheb2}
(x^2-1)U_{n-1}'(x)=nT_n(x)-xU_{n-1}(x), \qquad T_n^2(x)-(x^2-1)U^2_{n-1}(x)=1,
\end{equation}
after several simplifications, we obtain
$$
R_n'(x)=
%-n^2\frac{T_n(x)\left(T_{n}(x)-T_n(a)\right)-(x^2-1)U^2_{n-1}(x)}{\sqrt{1-x^2}\left(T_n(x)-T_n(a)\right)^2}=
-n^2\frac{1-T_n(a)T_{n}(x)}{\sqrt{1-x^2}\left(T_n(x)-T_n(a)\right)^2}.
$$
The denominator of the fraction does not vanish on $I_0$ (see $(c)$), consequently, $e_k$, $k=\overline{1,n}$, are solutions to the equality (\ref{equ_(b)}).
From properties of $T_n$ on $I_o$ and the inequality $0<1/T_n(a)<1$ for $a>1$, it follows that all $e_k$ belong to $I_o$. Now let us demonstrate that $e_k$ are the points $a_k$ forming the alternance. Alternation in sign of values $R_n(e_k)$ is deduced from (\ref{equ_(b)}). Indeed,  the polynomial $T_n(x)-1/T_n(a)$ attains positive and negative values alternately on $I$ due to properties of $T_n$. Furthermore,
$$
|R_n(e_k)|=n\frac{\sqrt{1-x^2}|U_{n-1}(e_k)|}{|T_{n}(e_k)-T_n(a)|}=
n\frac{\sqrt{1-T^2_n(e_k)}}{T_n(a)-T_{n}(e_k)}=
%n\frac{\sqrt{1-1/T^2_n(a)}}{T_n(a)-1/T_n(a)}=
\frac{n}{\sqrt{T^2_n(a)-1}}.
$$

The conclusion $(c)$ about the poles of the \textit{l.d.} (\ref{n.d.*}) follows from the property  of the \textit{l.d.} (\ref{n.d.Dan}) indicated before Theorem~\ref{th1} (see the text around the formula (\ref{ellipseChu})).
\end{Proof}
\begin{lemma}[Komarov \cite{KomarovPMA2013}]
\label{lemmaKom1}
Let
$$
P_n(x)=\frac{p'(x)}{p(x)},\qquad
Q_m(x)= \frac{q'(x)}{q(x)},\qquad m\le n,
$$
where  $p(x)=\prod_{k=1}^n (x-z_k)$ has only simple roots $z_k$ and $q(x)=\prod_{k=1}^m (x-\zeta_k)$. Then there are constants $\gamma_k$, such that
$$
P_n(x)-Q_m(x)\equiv \frac{p(x)}{q(x)}\sum_{k=1}^n \frac{\gamma_k}{(x-z_k)^2}.
$$
\end{lemma}

The following lemma was proved in \cite{KomarovPMA2013}  under the assumption of validity of a certain determinant identity. We show in \S\ref{s3} that the identity actually holds (Theorem~\ref{A}). Relying on it, we formulate the lemma in its final form.
\begin{lemma}
\label{lemmaKom2}
For $j,k=\overline{1,n}$, let $z_k$ be pairwise distinct real or complex conjugate numbers satisfying the condition
\begin{equation}
\label{usl_Kom}
|z_k|>1,
\end{equation}
and $c_j$ be pairwise distinct numbers belonging to $I$. Then the determinant of the matrix $A:=\left((c_j-z_k)^{-2}\right)$ does not vanish.
\end{lemma}
\begin{ProofTh1}
Let us assume that there exists a \textit{l.d.} $\rho^*_n(x):=\pi'(x)/\pi(x)$ of order $m\le n$ belonging to $\mathfrak{R}_n(a)$, such that $\|\rho^*_n\|^*<\|\tilde{\rho}^*_n\|^*$ on~$I$. We set $r_n(x):=\sqrt{1-x^2}\rho^*_n(x)$. According to Lemma~\ref{lemma4}, the function $R_n(x)=\sqrt{1-x^2}\tilde{\rho}^*_n(x)$ has an alternance on $I$ consisting of the points $a_k$, $k=\overline{1,n}$, in $I_o$,  therefore
$$
\textrm{sgn}(R_n(a_k)-r_n(a_k))=\pm (-1)^k, \qquad k=\overline{1,n}.
$$
Consequently, the continuous function $R_n-r_n$ (then the function $\tilde{\rho}^*_n-\rho^*_n$) reverses a sign at least $n$ times on $I_o$ and thus has at least $n-1$ pairwise distinct roots $c_j$, $j=\overline{1,n-1}$, in $I_o$.

Both $\tilde{\rho}^*_n$ and $\rho^*_n$ belong to $\mathfrak{R}_n(a)$, thereby there exist polynomials $p$ and $q$ of degree $n-1$ and $m-1$ correspondingly, such that for the denominators of $\tilde{\rho}^*_n$ and $\rho^*_n$ we have $T_n(x)-T_n(a)=p(x)(x-a)$ and $\pi (x)=q(x)(x-a)$. From this by Lemma~\ref{lemmaKom1} we get
$$
\tilde{\rho}_n^*(x)-\rho^*_n(x)=\left(\frac{p'(x)}{p(x)}+\frac{1}{x-a}\right)-\left(\frac{q'(x)}{q(x)}+\frac{1}{x-a}\right)
=\frac{p(x)}{q(x)}\sum_{k=1}^{n-1}\frac{\gamma_k}{(x-z_k)^2},
$$
where $z_k$ are the (simple) roots of $p$. Hence the points $c_j$, $j=\overline{1,n-1}$, are also zeros of the sum $\sum_{k=1}^{n-1}\gamma_k(x-z_k)^{-2}$ (the polynomial $p$ has no real roots on $I$ due to the conclusion $(c)$ of Lemma~\ref{lemma4}), and the system of linear equations
\begin{equation}
\label{syst1}
\sum_{k=1}^{n-1}\frac{\gamma_k}{(c_j-z_k)^2}=\tilde{\rho}_n^*(c_j)-\rho^*_n(c_j),\qquad j=\overline{1,n-1},
\end{equation}
with respect to the variables $\gamma_k$ is homogeneous and the determinants $\Delta_k$, $k=\overline{1,n-1}$ of the matrices with a replaced $k$th column are equal to zero. Let us recall that the poles of $\tilde{\rho}^*_n$ satisfy the condition~(\ref{usl_Kom}) here, owing to $(c)$ of Lemma~\ref{lemma4} and the assumption $a>\sqrt{2}$. Hence by Lemma~\ref{lemmaKom2} the determinant $\Delta$ of the matrix of (\ref{syst1}) does not vanish, consequently, the system~(\ref{syst1}) has the trivial solution $\gamma_k=\Delta_k/\Delta=0$, $k=\overline{1,n-1}$. Thus $\tilde{\rho}_n^*\equiv \rho^*_n$, which contradicts the assumption that there exists better approximation.

The latter conclusion of Theorem~\ref{th1} follows from $(b)$ of Lemma~\ref{lemma4}.
\end{ProofTh1}

\section{\textit{L.d.}s of least deviation from zero on $I$ in the norm $\|\cdot\|$}
\label{s2}

Several estimates for least deviation of \textit{l.d.}s from zero in the norm $\|\cdot\|$ were obtained in \cite{DDan2001}.
 In particular, it was shown there that if $\|\rho_n\|\le b^{-n-1}$ for some $b > 2$ then all poles of the \textit{l.d.} $\rho_n$ lie outside the ellipse $E_{p(A)}$ of the form (\ref{ellipseChu}) with $A=(b/2)^n$.
The following statement\footnote{As usually, for some positive sequences $\{a_n\}$ and $\{b_n\}$ \textit{the equivalence} $a_k\sim b_k$ means that $a_n/b_n\to 1$ for $n\to\infty$, and \textit{the weak equivalence} $a_n\asymp b_n$ does that there exist constants $\alpha$ and $\beta$ such that $0<\alpha\le a_n/b_n\le \beta<\infty$ for $n\ge n_0$.} is a supplement to this result.
\begin{theorem}
\label{th2}
For $n\ge 4$ and $a>\sqrt{2}\cdot(3\sqrt{n})^{1/n}$, let the \textit{l.d.} $\widetilde{\rho}_{n}(a;x)$ be least deviating from zero on $I$ in the norm $\|\cdot\|^*$ among all those belonging to $\mathfrak{R}_n(a)$. Then
$$
\|\widetilde{\rho}_{n}\| \asymp \frac{2n}{T_n(a)-T_{n-2}(a)}, \qquad n\ge 4.
$$
For $n\to \infty$ the sign $\asymp$ can be replaced by the sign $\sim$.
\end{theorem}

We first prove auxiliary results. The following lemma is an analogue of De la Vall\'{e}e Poussin's theorem (see \cite[Ch. 1, \S 2]{Dzyadyk}) for \textit{l.d.}s belonging to $\mathfrak{R}_n(a)$.
\begin{lemma}
\label{lemma3}
Let the \textit{l.d.} $\rho_n$ belong to $\mathfrak{R}_n(a)$ and have the poles $z_k$ being pairwise distinct and satisfying the condition~$(\ref{usl_Kom})$. Let there exist $n$ pairwise distinct points $t_k$ in $I$, such that  $\rho_n(t_k)=\pm(-1)^k \lambda_k$, $\lambda_k>0$, $k=\overline{1,n}$. Then for the \textit{l.d.} $\widetilde{\rho}_n$,  least deviating from zero on $I$ in the norm $\|\cdot\|$ among all those belonging to $\mathfrak{R}_n(a)$, 
$$
\min_{k=\overline{1,n}} \lambda_k \le \|\widetilde{\rho}_n\|\le \|\rho_n\|.
$$
\end{lemma}
\begin{Proof}
The right hand side inequality is obvious.  The non-existence of \textit{l.d.}s, belonging to $\mathfrak{R}_n(a)$, which norms are less than  $\min_{k=\overline{1,n}} \lambda_k$ (the left hand side inequality), follows by essentially the same method as in the proof of Theorem~\ref{th1}, and we do not produce it.
\end{Proof}
\begin{lemma}
\label{lemma5}
For $n\ge 4$, the \textit{l.d.}
\begin{equation}
\label{best_uniform}
\widetilde{\varrho}_{n}(a;x)=\frac{T_{n-1}(x)}{\int_{a}^x T_{n-1}(t)dt}, \qquad a>1+1/n,
\end{equation}
has the following properties.

$(a)$ For $a_k=\cos(\frac{k}{n-1}\pi)$, $k=\overline{0,n-1}$, in $I$, we have  $\widetilde{\varrho}_{n}(a_k)=\pm(-1)^k\lambda_k$, where $\lambda_k>0$ and
\begin{equation}
\label{estim_lambda}
\frac{2n}{T_{n}(a)-\frac{n}{n-2}T_{n-2}(a)+2\frac{n-1}{n-2}}
\le
\lambda_k
\le \frac{2n}{T_{n}(a)-\frac{n}{n-2}T_{n-2}(a)-2\frac{n-1}{n-2}}.
\end{equation}

$(b)$ The poles of $\widetilde{\varrho}_{n}$ lie in the closure of the ellipse $E_a$ of the form~$(\ref{ellipseChu})$. Moreover, if $t_n(a):=a\cdot(3\sqrt{n})^{-1/n}$ exceeds $1$ then all of them lie outside the closure of the ellipse $E_{t_n(a)}$ of the form $(\ref{ellipseChu})$.
\end{lemma}
\begin{Proof} We set
$$
Q_n(a;x):=\int\nolimits_{a}^x T_{n-1}(t)dt, \qquad f(x):=\tfrac{1}{n}T_{n}(x)-\tfrac{1}{n-2}T_{n-2}(x).
$$
The following well-known identity $Q_n(a;x)=\frac{1}{2}\left(f(x)-f(a)\right)$ is used below.
Throughout the proof $a>1+1/n$ and $n\ge 4$. Simple analysis of $f$ gives
\begin{equation}
\label{Q}
-\frac{n-1}{n(n-2)}\le Q_n(a;x)+\tfrac{1}{2}f(a)\le\frac{n-1}{n(n-2)}, \qquad x\in I.
\end{equation}
From this by the inequality $\frac{1}{2}f(a)>\frac{n-1}{n(n-2)}$ we get $Q_n(a;x)<0$ for $x\in I$.
Furthermore, in the point $a_k$, $k=\overline{0,n-1}$, being zeros of $\sqrt{1-x^2}U_{n-2}(x)$ (see~(\ref{Cheb1})), we have
$T_{n-1}(a_k)=\cos k\pi=(-1)^k$, $k=\overline{0,n-1}$ since $T_n(\cos \nu)=\cos n\nu$. It proves the alternation of signs of $\widetilde{\varrho}_{n}(a_k)$. The identity $|T_{n-1}(a_k)|=1$ and the inequality (\ref{Q}) yield the estimates (\ref{estim_lambda}).

Let us now prove $(b)$. We first show that all roots of $Q_n$ lie in the closure of the ellipse $E_a$ of the form~(\ref{ellipseChu}). To do so, we make sure that the inequality ${f(a)<|f(z)|}$ holds for all $z\in E_{a+\varepsilon}$, where $\varepsilon>0$. Taking into account the identity
\begin{equation}
\label{Cheb6}
T_n(z)= \tfrac{1}{2}\left((z+\sqrt{z^2-1})^n+(z-\sqrt{z^2-1})^n\right), \qquad z\in \mathbb{C},
 \end{equation}
where $\sqrt{z^2-1}$ is the branch, for which the value equals $1$ if ${z=\sqrt{2}}$, by direct calculations, we obtain that $|f(z)|$ attains a minimal value for $z=\pm r$ and a maximal value for $z=\pm i\sqrt{r^2-1}$ on the ellipse $\mathcal{E}_r:=\{z:|z+\sqrt{z^2-1}|=r\}$ with $r>1$ (i.e. when $z+\sqrt{z^2-1}=re^{i\varphi}$, $\varphi\in[0;2\pi)$). Since the ellipses
$E_{a+\varepsilon}$ and $\mathcal{E}_r$ coincide if $r=(a+\varepsilon)+\sqrt{(a+\varepsilon)^2-1}$, by (\ref{Cheb6}) and monotone increasing of $f(x)$ for $x>1$, we thus get
$$
|f(z)|_{z\in\mathcal{E}_r}\ge f(a+\varepsilon)>f(a),\qquad \varepsilon>0.
$$
Consequently, by Rouch\'{e}'s theorem, the polynomials $f(x)$ and $f(x)-f(a)$ have the same number of roots inside the ellipse $E_{a+\varepsilon}$. At the same time all roots of the polynomial $f$, as it can be easily shown,  belong to  $[-1-1/n;1+1/n]$, therefore to the closure of the ellipse $E_{a+\varepsilon}$ with $a>1+1/n$. By arbitrariness  of $\varepsilon$, all roots of $f(x)-f(a)$ (then all roots of $Q_n$) lie in the closure of $E_a$.

We now show that there are no roots of $Q_n$ in the closure of $E_{t_n(a)}$ of the form~(\ref{ellipseChu}) with $t_n(a)=a\cdot(3\sqrt{n})^{-1/n}$. It is only sufficient to prove that $|f(z)|<f(a)$ for $z\in E_{t_n(a)}$. To do so, we make the following observation. 
On the one hand, the extremal properties of $|f(z)|$ on $\mathcal{E}_r$, $r>1$, and the representation (\ref{Cheb6}) yield
$$
|f(z)|_{z\in E_{t_n(a)}}\le \tfrac{1}{n}T_n(t_n(a))+\tfrac{1}{n-2}T_{n-2}(t_n(a))<\tfrac{1}{n}(t_n(a)+\sqrt{t^2_n(a)-1})^n.
$$
From the other hand, $f(a)>(a+\sqrt{a^2-1})^{n}/(3n^{3/2})$, which is easy to check. These imply the strengthened inequality
$$
t_n(a)+\sqrt{t_n^2(a)-1}\le (3\sqrt{n})^{-1/n}(a+\sqrt{a^2-1}),
$$
that is valid for $t_n(a)$ mentioned, which, however, must exceed $1$ (it is so for ${a>(3\sqrt{n})^{1/n}>1+1/n}$). Taking into account that $f(a)$ is a constant, by Rouch\'{e}'s theorem, we finally conclude that the polynomial ${f(x)-f(a)}$ has no roots in the closure of $E_{t_n(a)}$.
\end{Proof}

\medskip

\begin{ProofTh2}
All the conclusions of Theorem~\ref{th2} follow from the analogue of De la Vall\'{e}e Poussin's theorem for \textit{l.d.}s (Lemma~\ref{lemma3}), where $\rho_n$ is the \textit{l.d.} $\widetilde{\varrho}_{n}$ of the form (\ref{best_uniform}), and  the estimates (\ref{estim_lambda}). The following observation should be only taken into account. Lemma~\ref{lemma3} requires the poles of $\widetilde{\varrho}_{n}$ to satisfy the condition~(\ref{usl_Kom}). It is easily seen from $(b)$ of Lemma~\ref{lemma5}
that the condition is satisfied if the  semi-minor axis  $\sqrt{t^2_n(a)-1}$ of the ellipse $E_{t_n(a)}$ exceeds $1$. It gives the condition $a> \sqrt{2}\cdot(3\sqrt{n})^{1/n}$, $n\ge 4$, in Theorem~\ref{th2}.
\end{ProofTh2}

\section{Corollary of main theorems}
\label{s-Col}
In approximation theory, Markov-Bernstein type estimates connecting values of polynomials and their derivatives are well-known (see for instance a wide survey in \cite[Ch. 5; Appendix A5]{BorwEld}). Now we give a result of this type, which is a corollary of Theorems~\ref{th1} and~\ref{th2}.
\begin{col}
Let $n\ge 4$ and $P_n(x)=(x-a)p(x)$, where $a>\sqrt{2}\cdot(3\sqrt{n})^{1/n}$ and $p=p(x)$ is a polynomial of degree $n-1$ having no roots in $I$. Then
\begin{equation}
\label{col}
\|P_n'\|^*\ge \frac{n\;\min_{x\in I}|P_n(x)|}{\sqrt{T_n^2(a)-1}}, \qquad
\|P_n'\|\ge \frac{2n \;\min_{x\in I}|P_n(x)|}{T_n(a)-T_{n-2}(a)+3}, \qquad n\ge 4.
\end{equation}

These inequalities are asymptotically precise in the sense that there exist polynomials
$\mathcal{P}_{n,1}$ and $\mathcal{P}_{n,2}$ such that for $n\to\infty$
$$
\frac{n\;\min_{x\in I}|\mathcal{P}_{n,1}(x)|}{\sqrt{T_n^2(a)-1}}\sim \|\mathcal{P}_{n,1}'\|^*,\qquad
\frac{2n\;\min_{x\in I}|\mathcal{P}_{n,2}(x)|}{T_n(a)-T_{n-2}(a)+3}\sim \|\mathcal{P}_{n,2}'\|.
$$
\end{col}
\begin{Proof} From Theorem~\ref{th1} under the above assumptions, it immediately follows that
$$
\frac{n}{\sqrt{T_n^2(a)-1}}\le \left\|\frac{P_n'}{P_n}\right\|^*\le
\frac{\left\|P_n'\right\|^*}{\min_{x\in I}|P_n(x)|}.
$$
Analogously Theorem~\ref{th2} (more precisely, the conclusion $(a)$ of Lemma~\ref{lemma5} after minor simplifications) yields the inequality for the norm $\|\cdot\|$. Note that the former inequality in (\ref{col}) holds already for $a>\sqrt{2}$ and $n\ge 1$.

We now prove the asymptotical precision. For $\mathcal{P}_{n,1}(x)=T_n(x)-T_n(a)$ we have $\left\|\mathcal{P}_{n,1}'\right\|^*=n\|U_{n-1}\|^*=n$, $\min_{x\in I}|\mathcal{P}_{n,1}(x)|=T_n(a)-1$, thereby
$$
\frac{n\;\min_{x\in I}|\mathcal{P}_{n,1}(x)|}{\sqrt{T_n^2(a)-1}\;\|\mathcal{P}_{n,1}'\|^*}=\sqrt{1-\tfrac{2}{T_n(a)+1}}= 1-o(1), \qquad o(1)>1,\qquad n\to\infty.
$$

Furthermore, for $\mathcal{P}_{n,2}(x)=\int_a^x T_{n-1}(t)\;dt$ it is true that
$\left\|\mathcal{P}_{n,2}'\right\|=\|T_{n-1}\|=1$,
$\min_{x\in I}|\mathcal{P}_{n,2}(x)|\ge \frac{1}{2n}(T_n(a)-\frac{n}{n-2}T_{n-2}(a)-3)$ for $n\ge 4$, and, consequently,
$$
\frac{2n\;\min_{x\in I}|\mathcal{P}_{n,2}(x)|}{\left(T_n(a)-T_{n-2}(a)+3\right)\,\|\mathcal{P}_{n,2}'\|}
\ge 1-\frac{2(\frac{1}{n-2}T_{n-2}(a)-3)}{T_n(a)-T_{n-2}(a)+3}= 1-o(1), \quad n\to\infty,
$$
where $o(1)>0$.\end{Proof}

\section{Criterion for best uniform approximation by \textit{l.d.}s}
\label{s3}

It is known that approximating properties of \textit{l.d.}s and polynomials are similar in many respects. For instance, on has for \textit{l.d.}s analogues of classical Jackson's, Bernstein's, Zygmund's, Dzyadyk's, and Walsh's theorems for polynomials \cite{Kosuh,DanDan2001,Novak}. However, there are fundamental differences as well. So, generally speaking, there exist no direct connection between an alternance and best approximation in the case of \textit{l.d.}s. Moreover, a \textit{l.d.} of best approximation can be non-unique. A pioneering example of this kind was given for $n=2$ in~\cite{DanKond2010}, then it was extended in \cite{KomarovNeed} to arbitrarily positive integer $n$. Connection between best approximation by \textit{l.d.}s and an alternance in the general case is still a question.  Nevertheless,  several analogues of Chebyshev's alternance theorem can be obtained for \textit{l.d.}s with poles under certain restrictions. In particular, the following statement is valid.

\begin{theorem}
\label{Crit}
Let all the poles of the \textit{l.d.} $\rho_n=\rho_n(x)$ be pairwise distinct and satisfy the condition $(\ref{usl_Kom})$.
In this case $\rho_n$ is a unique \textit{l.d.} of best uniform approximation of the continuous function $f=f(x)$ on $I$ if and only if there exist $n+1$ points forming an alternance of the difference $f-\rho_n$ on $I$.
\end{theorem}

This theorem was formulated and proved by Komarov in \cite{Komarov2014,KomarovPMA2013} under the essential assumption consisting in validity of the following determinant identity due to Borchardt (it was considered in \cite{Komarov2014,KomarovPMA2013} as a conjecture).

\begin{theorem}[Borchardt \cite{Borchardt,Mink}]
\label{A}
For $j,k=\overline{1,n}$, let  $\{z_k\}$ and $\{c_j\}$ be arbitrary disjoint collections of complex numbers and   $A:=\left((c_j-z_k)^{-2}\right)$ and ${B:=\left((c_j-z_k)^{-1}\right)}$ be matrices constructed on basis of them. Then
$$
\det A=\det B \cdot \; {\rm per\,} B,
$$
where ${\rm per\,} B$ is a permanent\footnote{
Let us recall the definition of a permanent. For any square matrix $M=(m_{i,j})$, $i,j=\overline{1,n}$, 
${{\rm per\,} M=\sum_{\sigma\in S_n} \prod_{i=1}^n m_{i,\sigma_i}}$, where $S_n=\{1,\ldots,n\}$  (cf. $\det M=\sum_{\sigma\in S_n}\textrm{sgn}(\sigma) \prod_{i=1}^n m_{i,\sigma_i}$).
} of the matrix $B$.
\end{theorem}

In conclusion we note that in the even earlier  paper \cite{Novak}, Novak, considering Theorem~\ref{A} also as a conjecture,  obtained an analogue of  Theorem~\ref{Crit} for the more narrow  class of \textit{l.d.}s with pairwise distinct real poles out of $I$.

\end{document}